\documentclass[10pt]{amsart}
\usepackage{amsmath,amssymb}

\newtheorem{thm}{Theorem}
\newtheorem{lem}[thm]{Lemma}
\newtheorem{prop}[thm]{Proposition}

\newcommand{\Z}{{\mathbb Z}}
\newcommand{\Q}{{\mathbb Q}}
\newcommand{\R}{{\mathbb R}}

\newcommand{\N}{{\mathbb N}}
\newcommand{\entropyA}{\eta_A}

\newcommand{\mat}{\operatorname{Mat}}
\newcommand{\ord}{\operatorname{ord}}

\newcommand{\OK}{{\mathcal O}_K}
\newcommand{\mnorm}{ {\mathbf N_{K/\Q}}}

\newcommand{\barQ}{{\bar{\Q}}}

\newcommand{\OS}{{\mathcal O}_S}

\begin{document}

\title{A lower bound for periods of matrices}

\author{ Pietro Corvaja, Z\'eev Rudnick and Umberto Zannier}
\address{Dipartimento di Mathematica e Inf., Via delle Scienze, 206,
  33100 Udine, Italy  \newline
({\tt corvaja@dimi.uniud.it})}
\address{Raymond and Beverly Sackler School of Mathematical Sciences,
Tel Aviv University, Tel Aviv 69978, Israel
({\tt rudnick@post.tau.ac.il})}
\address{ I.U.A.V. -DCA, S. Croce 191, 30135 Venezia, Italy 
({\tt zannier@iuav.it})}

\date{February 9, 2004} 

\begin{abstract}
For a nonsingular integer matrix $A$, we study the growth of the order
of $A$ modulo $N$. We say that a matrix is exceptional 
if it is diagonalizable, and a power  of the matrix
has all eigenvalues equal to powers of a single rational integer, or all
eigenvalues are powers of a single unit in a real quadratic field. 

For exceptional matrices, it is easily seen that 
there are arbitrarily large values of $N$
for which the order of $A$ modulo $N$ is logarithmically small. 
In contrast, we show that if the matrix is not exceptional, then the
order of $A$ modulo $N$ goes to infinity faster than any constant
multiple of $\log N$. 
\end{abstract}

\maketitle


\section{Introduction}

Let $A$ be a $d\times d$ nonsingular integer matrix, 
and $N\geq 1$ an integer.  
The order, or period, of $A$ modulo $N$ is defined as the least integer 
$k\geq 1$ such that $A^k=I \mod N$,  
where $I$ denotes the identity matrix. 
If  $A$ is not invertible modulo $N$ then we set 
$\ord(A,N) =\infty$ . 
In this note we study the minimal growth of $\ord(A,N)$ as $N\to
\infty$. 

If $A$  is of finite order (globally),  that is $A^r=I$ for some
$r\geq 1$, then clearly $\ord(A,N)\leq r$ is bounded.  
If $A$ is of infinite order, then  
$\ord(A,N)\to \infty$ as $N\to \infty$. 
Moreover,  in this case it is easy to see that
$\ord(A,N)$ grows at least logarithmically with $N$, in fact 
if no eigenvalue of $A$ is a root of unity 
then: 
$$
\ord(A,N)\geq \frac {d}{\entropyA} \log N +O(1)
$$
where $\entropyA:=\sum_{|\lambda_j|>1}\log |\lambda_j|$, the sum over all
eigenvalues $\{\lambda_j\}$ of $A$ which lie outside the unit circle 
($\entropyA$ is the entropy of the endomorphism of the
  torus $\R^d/\Z^d$ induced by $A$, or the logarithmic Mahler measure
  of the characteristic polynomial of $A$, and the condition that no
eigenvalue of $A$ is a root of unity is equivalent to ergodicity of
the toral endomorphism).

There are cases when the growth of $\ord(A,N)$ is indeed no faster
than logarithmic. 
For instance if we take $d=1$, and  $A=(a)$ 
 where $a> 1$ is an integer, and $N_k=a^k-1$ then 
 \begin{equation*}
\ord(A,N_k)=k  \sim  \frac {\log N_k}{\log a}   
 \end{equation*}
and so 
\begin{equation}\label{a mod N}
  \liminf \frac{\ord(A,N)}{\log N} = \frac 1{\log a}  < \infty
\end{equation}
in this case. 

The same behaviour occurs in the case of $2\times 2$ unimodular
matrices $A\in {\rm SL}_2(\Z)$ which are hyperbolic, that is $A$ has a pair
of distinct real eigenvalues $\lambda>1>\lambda^{-1}$. 
Then 
\begin{equation}\label{small order}
  \liminf \frac{\ord(A,N)}{\log N} = \frac {2}{\log \lambda} = \frac {2}{\entropyA} 
\end{equation}
See e.g. \cite{KR}
\footnote{ a special case of this appeared as a problem in the 54-th
W.L. Putnam Mathematical Competition, 1994, 
see  \cite[pages 82, 242]{Andreescu}).}.

These cases turn out to be subsumed by the following definition: 
We say that $A$ is {\em  exceptional} if it is of finite order or if 
it is diagonalizable and a power
$A^r$ of $A$ satisfies one of the
following:
\begin{enumerate}
{\em 
\item \label{cond a}
 The  eigenvalues of $A^r$ are all a power of 
a single rational integer $a>1$;
\item\label{cond b}
 The  eigenvalues of $A^r$ are all a power of 
a single unit $\lambda\neq \pm 1$ of a real quadratic field.
}
\end{enumerate}

We will see that if $A$ is exceptional, then 
there is some $c>0$ and  arbitrarily large 
integers $N$ for which $\ord(A,N)<c \log N$. 

Our main finding in this note is 
\begin{thm}\label{Thm1}
If $A\in \mat_d(\Z)$ is not exceptional   
then 
$$
\frac {\ord(A,N)}{\log N} \to \infty 
$$  
as $N\to \infty$. 
\end{thm}

A special case is that of diagonal matrices, e.g. 
$A=
\begin{pmatrix}
 a&0\\0&b 
\end{pmatrix}$. In that case Theorem~\ref{Thm1} says that 
$\ord(a,b;N)/\log N\to \infty$ if $a,b$ are multiplicatively
independent, in contrast with \eqref{a mod N}.

Theorem~\ref{Thm1} is in fact equivalent to a subexponential bound on
the greatest common divisor $\gcd(A^n-I)$ of the matrix entries of
$A^n-I$. We shall derive it from 
\begin{thm}\label{Thm2}
If $A\in \mat_d(\Z)$ is not exceptional 
then for all $\epsilon>0$
$$
\gcd(A^n-I) < \exp(\epsilon n)
$$
if $n$ is sufficiently large.
\end{thm}

In the  special case of a diagonal matrix such as $A=
\begin{pmatrix}
 a&0\\0&b 
\end{pmatrix}$, we have $\gcd(A^n-I) = \gcd(a^n-1,b^n-1)$.  
In \cite{BCZ} it is shown that if $a,b$ are
multiplicatively independent then for all $\epsilon>0$,
\begin{equation}\label{bcz} 
  \gcd(a^n-1,b^n-1) <\exp(\epsilon n)
\end{equation}
for $n$ sufficiently large, giving Theorem~\ref{Thm2} in that case.  
To prove Theorem~\ref{Thm2} in general, 
we will use a version of \eqref{bcz} for
$S$-units in a general number field \cite{CZ} .

We note that Theorem~\ref{Thm2} establishes upper bounds on
$\gcd(A^n-I)$. As for lower bounds, it is conjectured in \cite{AR}
that if $A$ has a pair of multiplicatively independent eigenvalues
then  $\liminf \gcd(A^n-I)<\infty$.

{\bf Motivation:} A natural object of study for number theorists, 
the periods of toral automorphisms  were also 
investigated by a number of  physicists and mathematicians 
interested  in  classical and quantum dynamics,  
see e.g. \cite{HB,  Keating, DF}. 
One reason for our own interest  also lies in the quantum dynamics 
of toral automorphisms: 
It has recently been shown that any 
ergodic automorphism $A\in {\rm SL}_2(\Z)$ of the $2$-torus admits
``quantum limits''  different from Lebesgue measure \cite{BFN}, if one
does not take into account the hidden symmetries (``Hecke operators'') 
found in \cite{KR-Hecke}. 
The key  behind the constructions of these measures
is the existence of values of $N$ satisfying \eqref{small order}, that is
$\ord(A,N)\sim 2\log N/\entropyA$. 
A higher-dimensional version of this would involve taking ergodic 
symplectic automorphisms $A\in {\rm Sp}_{2g}(\Z)$ of the $2g$-dimensional
torus. Theorem~\ref{Thm1} gives one obstruction  to extending 
the construction of \cite{BFN}  to the higher-dimensional case.

\section{Proof of Theorem~\ref{Thm2}}

Assume that for a certain positive
$\epsilon$
and all integers $n$ in  a certain infinite sequence ${\mathcal N}\subset\N$ we
have
\begin{equation}\label{eqno(1)} 
\gcd(A^n-I) >\exp(\epsilon n).
\end{equation}
We shall prove that $A$ is ``exceptional", in the sense of the above
definition.

We let $k\subset\overline\Q$ be the splitting field for the
characteristic polynomial of $A$, so we may put $A$ in Jordan form over
$k$, namely,
we may write
$$
A=PBP^{-1},
$$
where $P$ is an invertible $d\times d$ matrix over $k$ and $B$ is in Jordan
canonical form. \smallskip

For later reference we introduce a little notation related to the field $k$.

We let $M$ (resp. $M_0$) denote  the set of (resp. finite) places of $k$.
We shall normalize all the absolute values   {\it with respect to $k$},
i.e. in such
a way that the product formula $\prod_{\mu\in M}|x|_\mu=1$ holds for $x\in
k^*$, and
the absolute logarithmic Weil height reads  $h(x)=\sum_\mu
\log\max\{1,|x|_\mu\}$.
We also let $S$ be a finite
set of places of $k$ including the archimedean ones and we denote by
$\OS^*$ the
group of $S$-units in $k^*$, namely those elements $x\in k$ such that
$|x|_\mu=1$
for all $\mu\not\in S$.\smallskip

Note that $B^n-I=P^{-1}(A^n-I)P$; since the entries of $P$ and its inverse are
fixed independently of $n$, hence have bounded denominators as $n$ varies, this
formula shows that the entries of $B^n-I$ have a ``big'' g.c.d., in the
sense of
ideals of $k$, for $n\in\mathcal N$. Since the entries of $B^n-I$ are algebraic
integers, not necessarily rational, to express their g.c.d. we shall use the
formula-definition
$$
\log\gcd(B^n-I):=\sum_{\mu\in M_0}\log^-\max_{ij}|(B^n-I)_{ij}|_\mu,
$$
where $\log^-(x):=-\min(0,\log x)$; this is a nonincreasing nonnegative
function of $x>0$.\smallskip

Note that this definition agrees with the usual notion in case $B$ has rational
integer entries. From \eqref{eqno(1)} and the above formula
$B^n-I=P^{-1}(A^n-I)P$ we immediately
deduce that
\begin{equation}\label{eqno(2)}
\sum_{\mu\in M_0}\log^-\max_{ij}|(B^n-I)_{ij}|_\mu>
\frac{\epsilon}{ 2} n,
\qquad\hbox{for large $n\in\mathcal N$}\;.
\end{equation}
In fact, each entry of $B^n-I$ is a linear combination of entries of
$A^n-I$ with
coefficients having bounded denominators, whence
$|(B^n-I)_{ij}|_\mu\le c_\mu\max_{rs}|(A^n-I)_{rs}|_\mu$, where $c_\mu$ are
positive
numbers independent of $n$ such that $c_\mu=1$ for all but finitely many
$\mu\in M$. This proves \eqref{eqno(2)}.\smallskip

We start by showing that $B$ must be necessarily diagonal. In fact, if not some
block of $B$ would contain on the diagonal a $2\times 2$ matrix of the form
$$
\begin{pmatrix}\lambda & 1\\ 0&\lambda\end{pmatrix}
$$
where $\lambda$ is an (algebraic integer) eigenvalue  of $A$. Hence $B^n-I$
would
contain among its entries the numbers $\lambda^n-1$ and $\lambda^{n-1}n$. Then,
 for every $\mu\in M_0$, we would have
$$
\max_{ij}|(B^n-I)_{ij}|_\mu\ge \max (|\lambda^n-1|_\mu,|\lambda^{n-1}n|_\mu)
\ge |n|_\mu,
$$
whence $\log^-\max_{ij}|(B^n-I)_{ij}|_\mu\le \log^-|n|_\mu=-\log |n|_\mu$. In
conclusion,
$$
\sum_{\mu\in M_0}\log^-\max_{ij}|(B^n-I)_{ij}|_\mu\le
\sum_{\mu\in M_0}-\log |n|_\mu=\log n
$$
the last equality holding because of the product formula. However this
contradicts
\eqref{eqno(2)} for all  large $n\in\mathcal N$ and this contradiction proves that $B$ is
diagonal.\smallskip

Therefore from now on we assume that $B$ is a diagonal matrix formed with  the
eigenvalues $\lambda_1,\ldots ,\lambda_d$ of $A$, each counted with the
suitable multiplicity.

Another case now occurs when there exist two multiplicatively independent
eigenvalues, denoted $\alpha,\beta$. Now, from \eqref{eqno(2)} we get,
for large $n\in\mathcal N$,
\begin{equation} \label{eqno(3)} 
\sum_{\mu\in M_0}\log^-\max (|\alpha^n-1|_\mu,|\beta^n-1|_\mu)\ge
\sum_{\mu\in M_0}\log^-\max_{ij}|(B^n-I)_{ij}|_\mu>
\frac{\epsilon}{2}n \;.
\end{equation}
We are then in position to apply (after a little change of notation) the
following
fact from \cite{CZ}, stated as Proposition 2 therein:\smallskip

\begin{prop}[Proposition 2 of \cite{CZ}]
Let $\delta>0$. All but
finitely many  solutions $(u,v)\in(\OS^*)^2$
to the inequality
$$
\sum_{\mu\in M_0} \log^-\max\{|u-1|_\mu, |v-1|_\mu\}>\delta
\cdot \max\{h(u),h(v)\}
$$
satisfy one of finitely many  relations $u^av^b=1$, where $a,b\in\Z$ are
not both zero.
\end{prop}

Actually, Prop. 2 in [CZ] is a little stronger, since the summation is
over all $\mu\in M$ rather than the finite $\mu\in M_0$ and since it also
asserts that
the relevant pairs $(a,b)$ may be computed in terms of $\delta$. \smallskip

We apply this fact with $u=\alpha^n$, $v=\beta^n$ and $S$ containing the
finite set
of places of $k$ which are nontrivial on $\alpha$ or $\beta$; note
that \eqref{eqno(3)} 
implies
the inequality  of the proposition, with $\delta
=\epsilon/(2\max(h(\alpha),h(\beta)))$. We conclude that, for an infinity of
$n\in\mathcal N$, a same nontrivial relation $\alpha^{an}\beta^{bn}=1$ holds,
contradicting the multiplicative independence of $\alpha,\beta$.\smallskip

Therefore we are left with the case when all pairs of eigenvalues are
multiplicatively dependent. This means that they generate in $k^*$ a subgroup
$\Gamma$ of rank $\le 1$.

If the rank is zero  all the eigenvalues $\lambda_i$ are roots of unity, so the
matrix  $A$ has finite order and thus it is exceptional. Hence let us
assume from now
on that the rank is $1$. Let then $\lambda\in \Gamma$ be a generator
of the free part of $\Gamma$ (it exists by basic theory).  Then,  for
suitable
roots of unity $\zeta_1,\ldots,\zeta_d$ and rational
integers $a_1,\ldots ,a_d$ we may
write
\begin{equation}\label{eqno(4)}
\lambda_i=\zeta_i\lambda^{a_i},\qquad i=1,\ldots ,d.
\end{equation} 
Necessarily the $\zeta_i$ lie in $k$.

Let $\sigma$ be an automorphism of $k$. Then $\sigma$
fixes the set of eigenvalues, since $A$ is a   matrix defined over $\Q$; hence
$\sigma$ fixes the above group $\Gamma$. Let $r$ be the order of the torsion in
$\Gamma$, so the subgroup $[r]\Gamma$ of $r$-th powers in $\Gamma$ is cyclic,
generated by $\lambda^r$. (Note that  automatically $\zeta_i^r=1$ in
\eqref{eqno(4)}). Then 
$\sigma$   must send $\lambda^r$ to another generator of $[r]\Gamma$, whence
$$
\sigma(\lambda)^r=\lambda^{\pm r}.
$$
Therefore in particular $\lambda^r$ is at most quadratic  over $\Q$
(in fact, recall that $k/\Q$ is normal).\smallskip

Let us first assume that $\lambda^r$ is rational. Raising the
equations \eqref{eqno(4)} 
to the
power $2r$,   we see that the
eigenvalues $\lambda_i^{2r}$ of the matrix $A^{2r}$ are positive rationals;
since they
are algebraic integers, they are therefore positive rational integers.
Since they are
pairwise multiplicatively dependent they are powers of a same positive integer
(which can be taken $\lambda^{\pm 2r}$). We thus fall in another of the
exceptional  situations.\smallskip

The last case occurs when $\lambda^r$ is a quadratic irrational. Then some
automorphism $\sigma$ must send it to its inverse $\lambda^{-r}$. As
before, we may raise equations \eqref{eqno(4)} to the  $r$-th power to find
$\lambda_i^{r}=\lambda^{ra_i}$. Therefore
$\sigma(\lambda_i^{r})=\lambda_i^{-r}$. Since the $\lambda_i$ are algebraic
integers, the same is true for the    $\lambda_i^{\pm r}$, and hence we
find that
all the eigenvalues of $A^{r}$ are units (some of them possibly
equal to $\pm 1$) in a same quadratic field.

This concludes the proof.\medskip


\section{Proof of Theorem~\ref{Thm1}}
The following Lemma shows that Theorems~\ref{Thm1} and \ref{Thm2} are
in fact equivalent:   
\begin{lem}
Let $A$ be a nonsingular integer matrix of infinite order. 
Then the following are
equivalent:
\begin{enumerate}
  \item For all $\epsilon>0$, we have  $\gcd(A^n-I)< \exp(\epsilon n)$
  if $n$ is sufficiently large;
  \item  $\ord(A,N)/\log N\to\infty$.  
\end{enumerate}
\end{lem}
\begin{proof}
Assume that $\gcd(A^n-I)< \exp(\epsilon n)$ for all $\epsilon>0$. 
Fix $\epsilon>0$. 
Take $n=\ord(A,N)$ and note that  $N$ divides all the matrix entries
of $A^{\ord(A,N)}- I$. Since $A$ does not have finite order and thus 
$\ord(A,N)\to \infty$ as $N\to \infty$,
we have for $N$ sufficiently large that  
$$
N\leq  \gcd(A^{\ord(A,N)}- I ) < \exp(\epsilon \ord(A,N)) 
$$ 
Thus 
$$
\log N < \epsilon \ord(A,N)\;. 
$$
Since this holds for all $\epsilon>0$ we find $\ord(A,N)/\log N\to
\infty$.

Conversely,  suppose that there is some $\rho>0$ and an infinite
sequence of integers $\mathcal N$ so that $\gcd(A^n-I)>\exp(\rho n)$
for all $n\in \mathcal N$. Then for the sequence 
$N_n:=\gcd(A^n-I)$, $n\in \mathcal N$ (which is infinite since
$N_n>\exp(\rho n)$) we 
have 
$$
\ord(A,N_n)\leq n<\log\gcd(A^n-I)/\rho = \log N_n/\rho
$$ 
and thus  $\liminf \ord(A,N)/\log N <\infty$. 
\end{proof}

\section{Comments}

It is readily seen that  exceptional cases 
do in fact occur, and that they give rise to powers $A^n$ such that $\gcd
(A^n-I)$ is exponentially large, and hence to arbitrarily large
integers $N$ for which $\ord(A,N)$ is logarithmically small. 
The last case of the eigenvalues in a quadratic
field of course requires that the irrational ones occur in conjugate pairs,
since $A$ is defined over $\Q$, and that the determinant of $A$ is
$\pm 1$. Examples of such integer matrices can be produced from the
action of a fixed such $2\times 2$  hyperbolic matrix $A_0\in
SL_2(\Z)$ on tensor powers, or from $A_0\otimes \sigma$ where $\sigma$
is a permutation matrix. 

To see that the exceptional cases lead to exponentially large $\gcd$,
consider first the case that a power of $A$ has all eigenvalues a
power of a single integer $a>1$. As we have seen in the course of
proof of Theorem~\ref{Thm2}, replacing a matrix by a conjugate (over
$\barQ$) does not change the asymptotic behaviour. 
Thus we may assume that $A^r$ is diagonal with eigenvalues
$a^{m_1},\dots, a^{m_d}$. Then clearly $\ord(A^r,N)\leq \ord(a,N)$ and
taking $N_n:=a^n-1$ gives $\ord(a,N_n) = n \sim \log N_n/a$. Thus we
find $\ord(A,N_n) \leq r\log N_n/a$. 

Now assume that a power $A^r$ of $A$ has all its eigenvalues a power
of a single unit $\lambda>1$ in a real quadratic field $K$. Then 
for some matrix $P$ with entries in $K$, we have 
$A^r=PBP^{-1}$ with $B$ diagonal with eigenvalues $\lambda^{a_1},\dots
,\lambda^{a_d}$, where $a_i$ are integers which sum to zero.  

Since $P$ is only
determined up to a scalar multiple, we may, after multiplying $P$ by an
algebraic integer of $K$, assume that $P$ has entries in the ring of
integers $\OK$ of $K$, and then $P^{-1} =\frac 1{\det(P)} P^{ad}$
where $P^{ad}$ also has entries in $\OK$.

The entries of $ A^{rk}-I $ are thus $\OK$-linear combinations of   
$(\lambda^{a_i k}-1)/\det(P)$. 
We now note that
$$\lambda^{-k}-1 = -\lambda^{-k}(\lambda^k-1)$$
and thus the entries of $ A^{rk}-I $ are all $\OK$-linear combinations of
$(\lambda^{|a_i| k}-1)/\det(P)$, which are in turn $\OK$-multiples of
$(\lambda^k-1)/\det(P)$. In particular, $\gcd(A^{rk}-I)$, which is a
$\Z$-linear combination of the entries of $A^{rk}-I$, can be written as
$$
\gcd(A^{rk}-I) = \frac {\lambda^k-1}{\det(P)} \gamma_k
$$
with $\gamma_k\in \OK$.

Now taking norms from $K$ to $\Q$ we see
$$
|\gcd(A^{rk}-I) |^2  =
\frac{|\mnorm(\lambda^k-1)| }{|\mnorm(\det P)|}|\mnorm(\gamma_k)|
\;.
$$
Since $\gamma_k\neq 0$, we have $|\mnorm(\gamma_k)|\geq 1$ and thus
$$
|\gcd(A^{rk}-I) |^2 \geq \frac{ |\mnorm(\lambda^k-1)|} {|\mnorm(\det
 P)|}
\gg \lambda^k
$$
which gives $|\gcd(A^{rk}-I) | \gg \lambda^{k/2}$, namely exponential growth.

\end{document}